\RequirePackage{fix-cm}
\documentclass[twocolumn]{svjour3}          
\smartqed  
\usepackage{graphicx}

\usepackage{amsfonts,amsbsy,amscd,amsgen,dsfont,amsmath,amssymb,mathrsfs,amsthm}
\usepackage[colorlinks=true,allcolors=blue]{hyperref}
\usepackage{lipsum}
\usepackage{amsfonts}
\usepackage{graphicx}
\usepackage{epstopdf}
\usepackage{algorithmic}
\usepackage{authblk}

\usepackage{amsopn}

\usepackage{subfigure}
\usepackage{bm} 
\usepackage{xcolor,comment}
\usepackage{enumitem}
\usepackage[utf8]{inputenc}
\usepackage[T1]{fontenc}

\newcommand{\id}{\mathrm d}
\newcommand{\vc}{\mathbf}

\newcommand{\bxi}{\pmb{\xi}}

\renewcommand{\i}{\hat{i}}
\renewcommand{\tilde}{\widetilde}
\newcommand{\pard}[2]{\frac{\partial #1}{\partial #2}}

\DeclareMathAlphabet\mathbfcal{OMS}{cmsy}{b}{n}

\newcommand{\ds}{\displaystyle}
\newcommand{\R}{\mathbb{R}}

\newcommand{\C}{\mathbb{C}}
\newcommand{\Cn}{\C^n}

\newcommand{\grad}{\nabla} 
\newcommand{\Lag}{\mathcal L}

\newcommand{\Real}{\text{Re}}
\newcommand{\Imag}{\text{Im}}
\newcommand{\sech}{\text{sech}}

\theoremstyle{definition}
\newtheorem{thm}{Theorem}

\newtheorem{rem}{Remark}

\newtheorem{cor}{Corollary}

\newcommand{\cb}{\color{black}} 

\begin{document}

\title{Shape-morphing reduced-order models for nonlinear Schr\"odinger equations}

\titlerunning{Shape-morphing reduced-order models for NLS}

\author{William Anderson \and Mohammad Farazmand}

\authorrunning{W. Anderson \and M. Farazmand} 
\institute{W. Anderson \and M. Farazmand \at
Department of Mathematics, North Carolina State University, 2311 Stinson Drive, Raleigh, NC
27695-8205, USA \\ \email{wmander3@ncsu.edu, farazmand@ncsu.edu}
}

\date{Received: date / Accepted: date}

\maketitle

\begin{abstract}
We consider reduced-order modeling of nonlinear dispersive waves described by a class of nonlinear Schr\"odinger (NLS) equations. We compare two nonlinear reduced-order modeling methods: (i) The reduced Lagrangian approach which relies on the variational formulation of NLS and (ii) The recently developed method of reduced-order nonlinear solutions (RONS). First, we prove the surprising result that, although the two methods are seemingly quite different, they can be obtained from the real and imaginary parts of a single complex-valued master equation. Furthermore, for the NLS equation in a stationary frame, we show that the reduced Lagrangian method fails to predict the correct group velocity of the waves whereas RONS predicts the correct group velocity. Finally, for the modified NLS equation, where the reduced Lagrangian approach is inapplicable, the RONS reduced-order model accurately approximates the true solutions.

\end{abstract}

\section{Introduction}

The nonlinear Schr{\"o}dinger (NLS) equation models weakly dispersive nonlinear waves, and thus has been used to study both optical waves \cite{agrawal2013,berge1994} and surface water waves \cite{dysthe08,zakharov68}.
In particular, an interesting feature of NLS is its ability to reproduce rogue waves, i.e., waves of extreme magnitude which form spontaneously from a combination of random wave superposition and deterministic nonlinear interactions~\cite{dysthe08}. 
These rogue waves have been observed in both water waves~\cite{chabchoub2012a,chabchoub2011,dysthe08} and optical waves~\cite{akhmediev13,solli2007}. 

Because of adverse consequences of rogues waves, such as sinking ships or damaging off-shore platforms, their real-time prediction and statistical quantification are highly desirable~\cite{farazmand2019a}. 
As a result, several methods have been developed for rapid prediction of rogue waves using low-dimensional reduced-order models of NLS~\cite{Adcock12,adcock09,cousins15,cousins16,Desaix91,farazmand2017,PerezGarcia1996,ruban2015,ruban2015b}.

These reduced models all evolve a wave packet with time-varying length scale, amplitude, and phase. 
Different models are distinguished in how they evolve these variables over time. A prominent variational approach relies on the 
Lagrangian structure of NLS by substituting the wave packet into its infinite-dimensional Lagrangian to obtain a reduced finite-dimensional counterpart. 
Subsequently writing the classical Euler--Lagrange equations leads to a set of ordinary differential {\cb equations} (ODEs) for the wave packet variables.
We refer to this method as the \emph{reduced Lagrangian approach} and describe it in more detail in Section~\ref{sec:redLag} below. 
This variational approach was first developed in the optics community~\cite{Desaix91,PerezGarcia1996} and was later adopted for water waves~\cite{ruban2015,ruban2015b}.


More recently, Anderson and Farazmand~\cite{anderson2021} developed the reduced-order nonlinear solutions (RONS) which is applicable to a wide class of partial differential equations (PDEs), including NLS.
RONS obtains a set of reduced-order equations for the time-dependent parameters by minimizing the instantaneous error between the evolution of reduced-order model and the full dynamics of the underlying PDE. RONS does not require a Lagrangian formulation for the PDE and allows for enforcing conserved quantities of the PDE in the reduced-order model. We review this method in Section~\ref{sec:RONS} below.

The goals of this manuscript is to compare the performance of the reduced Lagrangian approach and RONS. We discover that, although the two methods follow very difference reduced-order modeling philosophies, they are intimately connected. More precisely, we prove that there exists a complex-valued master equation whose real part gives the RONS model whereas its imaginary part coincides with the reduced Lagrangian model. Even more surprisingly, for a large class of wave packets, the two models coincide. 

In spite of these intimate connections, there are notable differences between the two methods. In particular, we show that the reduced Lagrangian approach may fail to provide an ODE for certain 
wave packet variables. For instance, the reduced Lagrangian approach fails to predict the wave speed of traveling waves. On the other hand, RONS is guaranteed to produce well-posed reduced-order equations, and in particular, accurately predicts the speed of traveling waves. 

Finally, we consider the modified NLS (or MNLS) equation which is a higher-order approximation to water wave equations~\cite{dysthe79}. MNLS equation does not have a known Lagrangian formulation and as such the reduced Lagrangian approach is not applicable to it. However, RONS is still applicable to MNLS and produces accurate reduced-order models.

\subsection{Outline of the paper}
\begin{figure*}
	\centering
	\includegraphics[width=\textwidth]{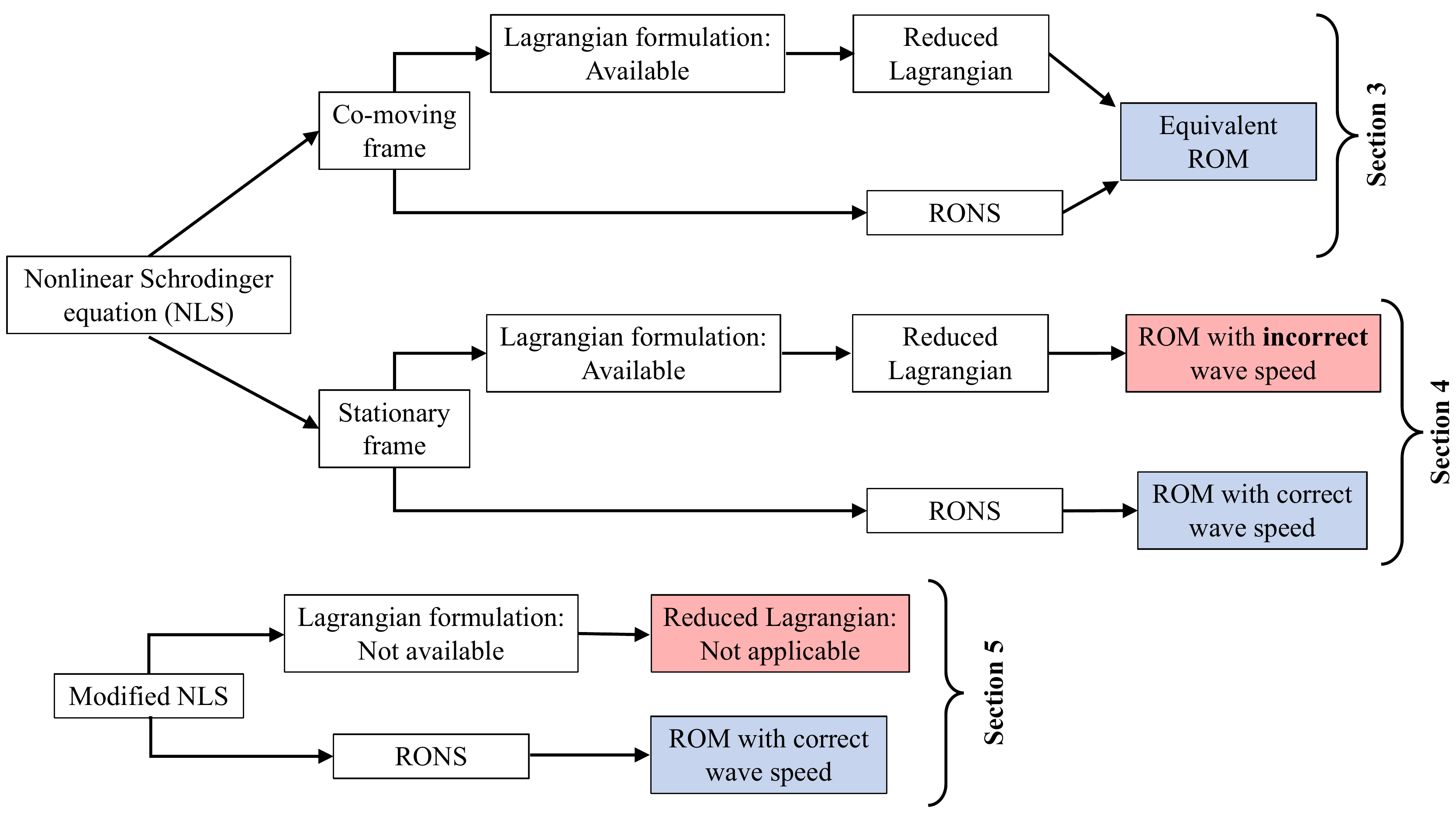}
	\caption{Summary of the main results. ROM stands for reduced-order model.}
	\label{fig:main_results}
\end{figure*}

This paper is organized as follows. In section \ref{sec:setup} we introduce the problem set-up and review the reduced Lagrangian approach and RONS. Section~\ref{sec:main_results} contains our main theoretical results where we show that RONS and the reduced Lagrangian approach are intimately connected. In section \ref{sec:lag_fails}, we show when the reduced Lagrangian approach fails while RONS produces reasonable results. In section \ref{sec:MNLS} we compare results produced by RONS with direct numerical simulation of MNLS. Section \ref{sec:conclusion} contains our concluding remarks. To facilitate the understanding of the main results, we summarize the outline of this paper in figure~\ref{fig:main_results}.

\section{Set-up and preliminaries}
\label{sec:setup}
Nonlinear Schr\"odinger equation is a perturbative model which describes the propagation of waves in both optics and hydrodynamics~\cite{zakharov68}.
In one {\cb dimension}, slowly varying modulations $\tilde u(x,t)$ to a carrier wave $\exp[{\cb i}(k_0x-\omega_0t)]$,
satisfy the NLS equation
\begin{equation}
\frac{ \partial u }{ \partial t  } = -\frac{1}{2} \frac{\partial u}{ \partial x } - \frac{ i }{ 8 } \frac{ \partial^2 u }{ \partial x^2  } - \frac{ i }{ 2 } |u|^2u,
\label{eq:NLS_advection}
\end{equation}
where we have non-dimensionalized the equation using the rescaled variables  $t = \omega_0\tilde{t}$, $x = k_0\tilde{x}$, and $u = k_0\tilde{u}$. 
The variables $\tilde{t}, \ \tilde{x}$, and  $\tilde{u}$ respectively represent time, space, and the wave envelope. The constants $k_0$ and $\omega_0$ denote the wave number and angular frequency of the carrier wave, respectively. 
It is common to consider a co-moving frame, translating with the constant speed $1/2$, so that the advection term is eliminated and NLS becomes
\begin{equation}
\frac{ \partial u }{ \partial t  } = - \frac{ i }{ 8 } \frac{ \partial^2 u }{ \partial x^2  } - \frac{ i }{ 2 } |u|^2u.
\label{eq:NLS}
\end{equation}

The NLS equation can also be derived from a variational formulation, where the {\cb extremizers} of the action functional $\int \mathcal L\id t$
satisfy the NLS equation~\cite{sulem2007}. The Lagrangian is given by
\begin{align}
\Lag = \int_{\mathbb R} \bigg( &\frac{1}{2} \big( iu_tu^\ast-iuu_t^\ast \big) + \frac{\gamma}{4} \big( iu_xu^\ast-iuu_x^\ast \big) \nonumber \\
 &+ \frac{1}{8} |u_x|^2 - \frac{ 1  }{ 4 } |u|^{4} \bigg)  \id x.
\label{eq:NLS_Lag_general}
\end{align}
This Lagrangian depends on a parameter $\gamma\in\{0,1\}$.
When $\gamma =1$, the corresponding Euler-Lagrange equations coincide with the NLS equation~\eqref{eq:NLS_advection} in the stationary frame.
Similarly, $\gamma = 0$ leads to the NLS equation~\eqref{eq:NLS} in the co-moving frame.

Now we turn our attention to reduced-order models that seek to approximate the solutions of NLS. In particular, we consider two classes of methods. The first class consists of methods that derive the reduced model from the Lagrangian~\eqref{eq:Lag_general}.
We review this method in Section~\ref{sec:redLag}. The second method is based on reduced-order nonlinear solutions (RONS) which was 
developed recently for reduced-order modeling of general PDEs~\cite{anderson2021}. RONS does not require a Lagrangian formulation and allows for enforcing conserved quantities of the PDE in the reduced model. 
We review RONS in Section~\ref{sec:RONS}

\subsection{Reduced Lagrangian approach}\label{sec:redLag}
In this section, we describe reduced-order models that take advantage of the known Lagrangian for NLS~\cite{Desaix91,PerezGarcia1996,ruban2015,ruban2015b}.
The first step in this reduced Lagrangian approach is to choose a suitable ansatz $\hat u (x, \vc q(t))$ for a wave packet with prescribed 
dependence on the spatial variable $x$ and a set of time-dependent parameters $\vc q(t)$. These parameters determine the
attributes of the wave packet such as its amplitude, wavelength and phase.

The next step is to derive a set of equations for
evolving the variables $\vc q(t)$ so that the resulting reduced-order solution $\hat u (x, \vc q(t))$ approximates a true solution of NLS.
A commonly used ansatz is the Gaussian wave packet,
\begin{equation}
\hat u(x, \vc q(t)) = A(t)\exp\bigg[- \frac{x^2}{L^2(t)} + i\frac{ x^2 U(t) }{ L(t) } + i \phi(t) \bigg],
\label{eq:gauss_ansatz}
\end{equation}
which has previously been considered by Desaix et al. \cite{Desaix91}, P{\'e}rez-Garc{\'i}a et al. \cite{PerezGarcia1996}, and Ruban \cite{ruban2015,ruban2015b}. Here $\vc q = (A, L, U, \phi)$ comprises the time-dependent variables, where $A$ is the wave amplitude, $L$ is a length scale which controls the wave width, $U$ determines the phase modulation, and $\phi$ is the phase.

To derive the evolution equations for the variables $\vc q(t)$, one first substitutes the ansatz~\eqref{eq:gauss_ansatz} into the Lagrangian~\eqref{eq:NLS_Lag_general} to obtain the \emph{reduced Lagrangian},
\begin{align}
	\hat \Lag (\vc q, \dot{  \vc q } ) = &\frac{\sqrt{\pi } A^2 }{16L} \bigg[ \sqrt{2} \Big( -2 L^3 \dot{U}  \nonumber \\ 
	&+ L^2 \Big(2 \dot{L} U+U^2-8 \dot{\phi } \Big) +1\Big) -2 A^2 L^2 \bigg].
	\label{eq:reduced_lag_comoving}
\end{align}

Note that the reduced Lagrangian $\hat{\mathcal L}$ depends only on the parameters $\vc q(t)$ and their time derivative $\dot{\vc q}(t)$ since the integration over space eliminates the spatial variable $x$. 
Writing the classical Euler-Lagrange equations for the reduced Lagrangian~\eqref{eq:reduced_lag_comoving}, 
we obtain
	\begin{align}
	&\dot{A} = \frac{A U}{4 L},&  &\dot{L} = -\frac{U}{2},& \nonumber \\
	&\dot{U} = \frac{\sqrt{2} A^2 L^2-2}{4 L^3},&  \quad &\dot{\phi} = \frac{1}{4 L^2}-\frac{5 A^2}{8 \sqrt{2}}.&
	\label{eq:param_odes}
	\end{align}
These ODEs provide a set of equations for evolving the parameters of the ansatz~\eqref{eq:gauss_ansatz}.

In the above derivation, we have used the Lagrangian in the co-moving frame ($\gamma=0$), as is customary. The reduced Lagrangian approach can also be applied to the problem in the stationary frame ($\gamma=1$), as we discuss in Section~\ref{sec:lag_fails}.

\begin{rem}\label{rem:NLS_E}
	There is no guarantee that the reduced Lagrangian approach will produce equations for all ansatz parameters. 
	As an example, we will show that by making a slight modification to the Gaussian ansatz the reduced Lagrangian approach fails to produce an equation for the phase $\phi$.
	
	Since NLS conserves mass,
     $ \int_{\R} | u(x,t) |^2  \id x$,
	it is desirable to enforce that the reduced solution $\hat u$ also conserves this quantity.
	As a result, this conservation law is often enforced explicitly by requiring 
	$$A(t) = A(0)\sqrt{L(0)/L(t)}$$ 
	so that $A$ is no longer an independent parameter in the ansatz.
	However, rewriting $A(t)$ as a function of $L(t)$ in this manner leads to a reduced Lagrangian which fails to produce any equations for the evolution of $\phi$ (see, e.g., Ruban~\cite{ruban2015,ruban2015b}).
\end{rem}

\subsection{RONS reduced-order model}\label{sec:RONS}
\begin{figure*}
	\centering
	\includegraphics[width=0.9\textwidth]{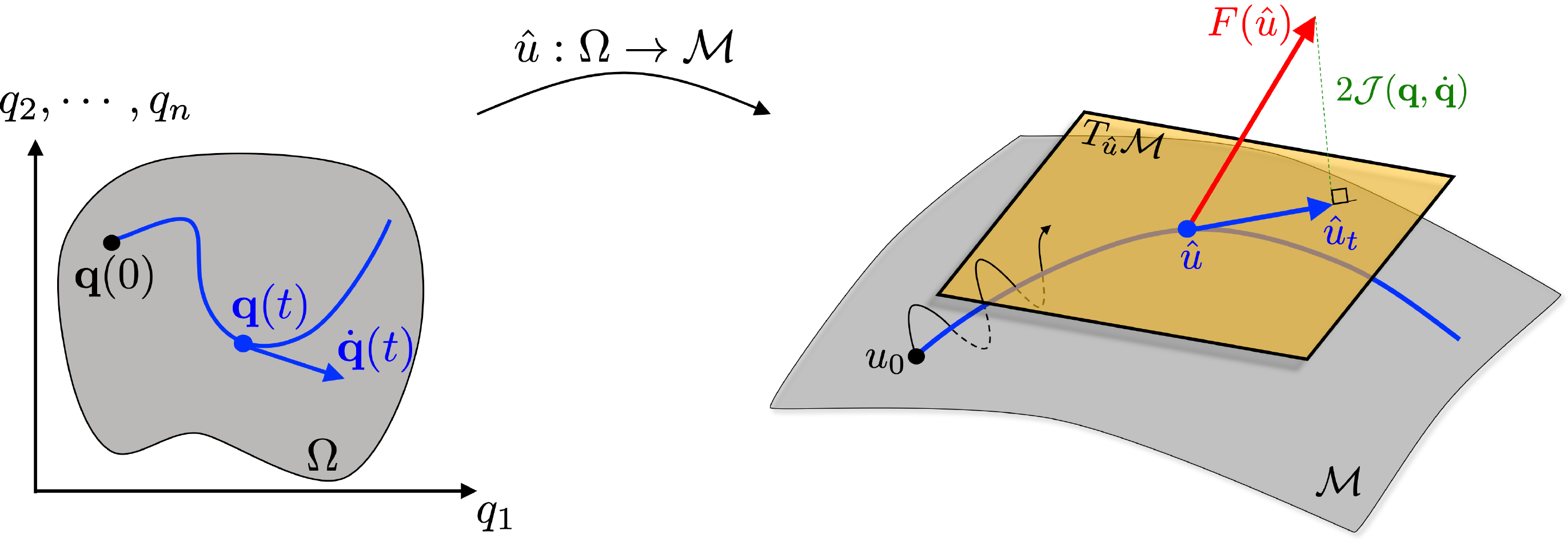}
	\caption{Geometric illustration of RONS. The ansatz $\hat u$ maps the parameters $\vc q$ to the ansatz manifold $\mathcal M$.}
	\label{fig:Geometric_Explanation}
\end{figure*}

In this section, we briefly review the method of reduced-order nonlinear solutions (RONS). A more detailed description of the method can be found in \cite{anderson2021}.

RONS is applicable to general PDEs of the form 
\begin{equation}
\frac{ \partial u }{ \partial t } = F(u), \quad u( \vc x, 0 ) = u_0(\vc x),
\label{eq:general_pde}
\end{equation}
where $u:D \times \R^+ \rightarrow \C$ is the solution, $D \subseteq \R^p$ denotes the spatial domain,
$F$ is a nonlinear differential operator, and $u_0$ is the initial condition.
For the case of one dimensional NLS, we have $p = 1$ and $D=\mathbb R$.
We assume $u(\cdot,t)$ belongs to a Hilbert space $H$ with the inner product $\langle \cdot,\cdot \rangle_H$
and the corresponding norm $\|\cdot \|_H$. For all RONS results in this manuscript we take the Hilbert space $H$ to be the space of complex square-integrable functions. 

Similar to the reduced Lagrangian case, RONS requires an ansatz $\hat u(x,\vc q (t) )$ for the shape of 
the reduced-order solution. Once an ansatz is chosen, RONS produces a set of equations for evolving the
time-dependent variables $\vc q(t)$. The RONS equations are obtained by minimizing the discrepancy between the reduced-order dynamics and the true PDE dynamics.

More precisely, we define the instantaneous error
\begin{equation}
\mathcal J( \vc q, \dot{ \vc q }) = \frac{1}{2} \|\hat u_t-F(\hat u)\|_H^2,
\label{eq:costfunctional}
\end{equation}
which measures the discrepancy between the rate of change of the ansatz $\hat u_t$ and the true rate dictated by the right-hand side of the PDE $F(\hat u)$.
To obtain the optimal evolution of the parameters $\vc q(t)$, we solve the minimization problem
\begin{equation}\label{eq:constJ}
	\min_{\dot{\vc q} \in \R^n} \mathcal J(\vc q,\dot{\vc q}),
\end{equation}
which minimizes the error $\mathcal J$ over all possible $\dot{\vc q}$ for a given $\vc q$.
The unique solution to this optimization problem is given by
\begin{equation}
\dot{\vc q} = [ M^{ (r) }(\vc q) ]^{-1}   \vc f^{ (r) } (\vc q),
\label{eq:qdot_const}
\end{equation}
which we refere to as the RONS equation.
Here $M\in \C^{n\times n}$ is the Hermitian positive definite \textit{metric tensor} whose entries are given by
\begin{equation}
M_{ij} = \left\langle \pard{\hat u}{q_i}, \pard{\hat u}{q_j}\right\rangle_H,\quad i,j\in\{1,2,\cdots,n\},
\label{eq:metricT}
\end{equation}
We denote the real and imaginary parts of the metric tensor by $M^{(r)}$ and $M^{(i)}$, respectively, so that
$M = M^{ (r) } + i M^{ (i) }$.
The vector field $\vc f:\Cn\to \Cn$ is defined by
\begin{equation}
f_i = \left\langle \pard{\hat u}{q_i}, F(\hat u)\right\rangle_H,\quad i=1,2,\cdots, n.
\end{equation}
We write this vector field in terms of its real and imaginary parts, $\vc f = \vc f^{(r)}+i\vc f^{(i)}$. 
Note that the RONS equation~\eqref{eq:qdot_const} involves only the real part of the the metric tensor $M$ and
the vector field $\vc f$.


A geometric illustration of this minimization process is shown in figure \ref{fig:Geometric_Explanation}. 
The ansatz has a fixed spatial dependence and therefore we can view $\hat u(\cdot,\vc q)$ as a map from the parameters $\vc q$ to the Hilbert space $H$.
Thus, $\hat u$ maps the set of all viable parameter values $\Omega\subseteq\mathbb R^n$  to an $n$-dimensional manifold $\mathcal M$ which lies in the Hilbert space.
An arbitrary and smooth parameter evolution defines a curve  $\vc q(t)$ where, at any point on the curve, its tangent vector is $\dot{  \vc q }(t)$. 
The curve representing an evolution of parameters $\vc q(t)$ is mapped to a curve which lies on $\mathcal M$, and the tangent vector $\dot{  \vc q }(t)$ is mapped to a tangent vector $\hat u_t$ on the tangent space $T_{\hat u} \mathcal{M}$ of the manifold. 
Since the ansatz manifold $\mathcal M$ is not necessarily invariant under the dynamics of the PDE~\eqref{eq:general_pde}, $F(\hat u)$ will not generally belong to the tangent space $T_{\hat u} \mathcal{M}$. 
Solving the minimization problem~\eqref{eq:constJ} is equivalent to orthogonal projection of $F(\hat u )$ onto $T_{\hat u} \mathcal{M}$.

As we show in~\cite{anderson2021}, and review in section~\ref{sec:conservedQuan} below, if the PDE has certain conserved quantities, the RONS equation can be readily modified to ensure that the reduced-order model respects these conserved quantities.
Furthermore, for the special case of a linear ansatz,
\begin{equation}
\hat u(\vc x,\vc q (t)) = \sum_{i = 1}^n q_i(t) u_i(\vc x),
\label{eq:lin_ans}
\end{equation}
where $\{u_i\}_{i=1}^{n}$ are a set of time-independent modes, the RONS equations coincide with standard Galerkin truncation. In this special case, the ansatz manifold $\mathcal M$ becomes a linear subspace of the Hilbert space $H$ and the Galerkin truncation becomes the orthogonal projection onto this linear subspace.
We refer to Ref.~\cite{anderson2021} for a detailed explanation.

{\cb As a special case, the modes $u_i$ in the linear ansatz~\eqref{eq:lin_ans} can be obtained through proper orthogonal decomposition (POD), which is described for NLS in \cite{shlizerman2011}. When using a linear ansatz with modes obtained from POD, the reduced-order equations produced by RONS coincide with those of classical POD-Galerkin reduced-order models. However, Galerkin truncation is not applicable to the nonlinear Gaussian ansatz~\eqref{eq:gauss_ansatz}.
As shown in \cite{anderson2021}, RONS with the Gaussian ansatz provides a much more accurate reduced-order solution to NLS as compared to classical POD-Galerkin techniques with comparable number of modes (see figure 4.2 in~\cite{anderson2021}).
}

Applying RONS to the NLS equation~\eqref{eq:NLS} in the co-moving frame using the Gaussian ansatz~\eqref{eq:gauss_ansatz}, we obtain a set of ODEs
\begin{align}
&\dot{A} = \frac{A U}{4 L},&  &\dot{L} = -\frac{U}{2},& \nonumber \\
&\dot{U} = \frac{\sqrt{2} A^2 L^2-2}{4 L^3},&  \quad &\dot{\phi} = \frac{1}{4 L^2}-\frac{5 A^2}{8 \sqrt{2}},&
\label{eq:RONS_NLS_eqn}
\end{align}
for the evolution of the parameters $\vc q = (A,L,U,\phi)$.
These ODEs coincide exactly with equation~\eqref{eq:param_odes} produced by the reduced Lagrangian approach. This is quite remarkable and unexpected since RONS and reduced Lagrangian constitute two distinct approaches to reduced-order modeling. In particular, RONS does not use the Lagrangian structure of NLS at all. Nonetheless, these two distinct approaches lead to the same set of reduced-order equations.

It is natural therefore to ask whether this is just a coincidence or there is a deeper connection between 
RONS and the reduced Lagrangian approach. In Section~\ref{sec:main_results}, we show that in fact there is an unexpectedly intimate connection between RONS and the reduced Lagrangian approach. 
Furthermore, although the reduced Lagrangian approach and RONS do not always produce the same reduced-order equations, we prove that the two methods coincide for a broad class of ansatz $\hat u$.

{\cb
\subsection{Enforcing conserved quantities}\label{sec:conservedQuan}
Some PDEs have conserved quantities which are invariant along their trajectories, and it is highly desirable that reduced-order models preserve these quantities as well. 
Invariant quantities of the PDE represent physical properties of true solutions which should also be present in the reduced-order solutions, and failing to enforce these conserved quantities may lead to reduced-order models which produce nonphysical results \cite{majda2012}. This fact has renewed interest in reduced-order models that preserve conserved quantities; see, for instance,  Ref. \cite{peng16} for Hamiltonian systems and Ref. \cite{karasozen18} for NLS.

The reduced Lagrangian approach does not provide a systematic method for ensuring preservation of conserved quantities in the reduced-order equations. 
However, in RONS, one can easily enforce any number of conserved quantities. 
A detailed description of this procedure is given in \cite{anderson2021}. 
Here we briefly summarize how to enforce conserved quantities within the RONS framework and show that the reduced-order equations for NLS~\eqref{eq:param_odes} preserve mass and energy.

Consider $m$ conserved quantities of PDE~\eqref{eq:general_pde}, $I_k: H \to \R$ with $k\in\{1,2,\cdots,m\}$. 
These quantities must remain invariant for all time so that $I_k(u( \cdot ,t)) = I_k(u_0)$.
 To obtain parameter ODEs which preserve conserved quantities, we solve the constrained minimization problem
\begin{align}\label{eq:min_cons}
& \min_{\dot{\vc q} \in \R^n} \mathcal J(\vc q,\dot{\vc q}),\nonumber\\
&\mbox{s.t.}\quad I_{ k }(\vc q (t))= I_k (\vc q(0))  , \quad k =1,2,...,m, \quad \forall t\geq 0,
\end{align}
where we write $I_{k}( \vc q(t) )$ as shorthand for $I_k(\hat u( \cdot ,\vc q (t) ) )$.
If a solution to the contrained minimization problem~\eqref{eq:min_cons} exists, it must satisfy
\begin{equation}\label{eq:min_cons_soln}
\dot{\vc q} = [ M^{ (r) }(\vc q) ]^{-1} \left[  \vc f^{ (r) } (\vc q) -\sum_{ k = 1 }^{ m } \lambda_k \nabla I_k(\vc q)\right].
\end{equation}
Here $M$ and $\vc f$ are the same metric tensor and vector field which appear in the solution to the unconstrained problem~\eqref{eq:qdot_const}. The Lagrange multiplier $\pmb \lambda=(\lambda_1,\lambda_2,\cdots,\lambda_m)^\top$ is the unique solution to the linear equation
\begin{equation}
C \pmb \lambda = \vc b,
\label{eq:lambda_eq}
\end{equation}
where $C(\vc q)\in \R^{m\times m}$ is the symmetric positive definite \textit{constraint matrix}
\begin{equation}
C_{ij} = \langle \grad I_j ,[M ^{ (r) }]^{-1} \grad I_i\rangle,\quad i,j\in\{1,2,\cdots,m\}, 
\label{eq:constMat}
\end{equation}
and the vector $\vc b( \vc q ) = (b_1,b_2,\cdots,b_m)^\top\in\R^m$ has entries defined by
\begin{equation}
b_i = \langle \grad I_i ,[M^{ (r) }]^{-1} \vc f^{ (r) } \rangle.
\label{eq:b}
\end{equation}
The inner product $\langle \cdot, \cdot \rangle$ is the standard Euclidean inner product.

Here we consider two important conserved quantities of NLS
\begin{align}
I_1(u) &= \int_{ \R }  |u|^2 \ \id x, \nonumber \\
I_2(u) &= \frac{ 1 }{ 8 }\int_{ \R }|u_x|^2 \ \id x - \frac{ 1 }{ 4 } \int_{ \R } |u|^4 \ \id x,
\end{align}
which represent mass and energy of the system, respectively.
These conserved quantities evaluated at the Gaussian ansatz~\eqref{eq:gauss_ansatz} become
\begin{align}
I_{1}( \vc q )&= \sqrt{\frac{\pi }{2}} A^2 L,  \nonumber\\
I_{2}( \vc q ) &= \frac{\sqrt{\pi } A^2 \left(\sqrt{2} \left(L^2 U^2+1\right)-2 A^2 L^2\right)}{16 L} .
\end{align}
Interestingly, the corresponding equation~\eqref{eq:b} yields $\vc b = (0,0)^\top$. Since the constraint matrix $C$ is symmetric positive definite, the only value of $\pmb \lambda$ which satisfies equation~\eqref{eq:lambda_eq} is $\pmb \lambda =  (0,0)^\top$.
So, in this case the RONS reduced-order equations for the unconstrained problem~\eqref{eq:constJ} and constrained problem~\eqref{eq:min_cons} coincide, meaning that total mass and energy are automatically preserved by the reduced-order equations~\eqref{eq:param_odes} for NLS.
}

\section{Co-moving frame: reduced-Lagrangian and RONS models coincide}\label{sec:main_results}

In this section, we show that the reduced-order equations produced by the reduced Lagrangian approach and RONS can in fact be derived from the same nonlinear system of equations. 
More specifically, there is a complex-valued master equation whose real part yields the RONS reduced-order model and its imaginary part leads to the reduced Lagrangian model. 
We then prove that, for a broad class of ansatz $\hat u$, the two models are identical, i.e., the real and imaginary parts of the master equation coincide.

In fact, our results can be proved for a class of PDEs slightly more general than NLS. 
To state the results for the most general case, consider the Lagrangian
\begin{equation}
	\Lag = \int_{\mathbb R} \bigg( \frac{1}{2} \big( iu_tu^\ast-iuu_t^\ast \big) + h(u, u_x, u^*,  u_x^*) \bigg)  \id x,
	\label{eq:Lag_general}
\end{equation}
where the map $h$ satisfies
\begin{align}
	\bigg( \frac{ \partial h }{ \partial u } \bigg)^* = \frac{ \partial h }{ \partial u^* }. 
\end{align}
Taking the first variation of the corresponding action functional leads to the PDE
\begin{equation}
	\frac{\partial u}{\partial t} = -i \frac{\partial }{ \partial x } \frac{\partial  h }{ \partial u_x^* } + i \frac{\partial h }{ \partial u^* }
	\label{eq:PDE_gen}
\end{equation}
Note that for the special case where
\begin{equation}
	h = \frac{\gamma}{4} \big( iu_xu^\ast-iuu_x^\ast \big)+\frac{1}{8}|u_x|^2 - \frac{1}{4} |u|^4, 
\end{equation}
we retrieve the Lagrangian~\eqref{eq:NLS_Lag_general}
and the corresponding PDE coincides with the NLS equation.

The following theorem reveals the close relationship between the reduced-order equations produced by RONS and the reduced Lagrangian approach.
\begin{thm}\label{thm:nonlin_sys}
	Consider PDE~\eqref{eq:PDE_gen} and
	let $\hat u (x, \vc q(t))$  denote a sufficiently smooth ansatz. Define the master equation
	\begin{equation}
	M ( \vc q ) \dot{ \vc q } =  \pmb \xi ( \vc q ) + \pmb \eta ( \vc q ) ,
	\label{eq:nonlin_sys}
	\end{equation} 
	where $M$ is the metric tensor~\eqref{eq:metricT} and $\pmb \xi:\Cn\to \Cn, \ \pmb \eta:\Cn\to \Cn$ are vector fields whose entries are
	\begin{equation*}
	\begin{aligned}
		\xi_k &= -i \int_{\R} \frac{\partial \hat u}{ \partial q_k}  \frac{\partial h}{ \partial \hat u }  \ \id x,\\
		  \eta_k &= -i \int_{\R} \frac{\partial \hat u_x}{ \partial q_k} \frac{\partial h}{ \partial \hat u_x } \ \id x,\quad k=1,2,\cdots, n.
	\end{aligned}
	\end{equation*} 
Then the following hold:
\begin{enumerate}
	\item The real part of~\eqref{eq:nonlin_sys} coincides with the reduced-order equations produced by RONS.
	\item The imaginary part of~\eqref{eq:nonlin_sys} coincides with the reduced-order equations produced by the reduced Lagrangian approach.
\end{enumerate}
	 
\end{thm}
\begin{proof}
	See Appendix \ref{sec:nonlin_sys_proof}.
\end{proof}

We emphasize that the NLS equations~\eqref{eq:NLS_advection} and~\eqref{eq:NLS} satisfies the assumptions of Theorem \ref{thm:nonlin_sys} and so this result is applicable when comparing RONS and the reduced-Lagrangian approach for NLS.

\begin{rem}\label{rem:HPD}
	The metric tensor $M$ which appears in equation~\eqref{eq:nonlin_sys} is Hermitian positive-definite. Thus, $\Real[M]$ is symmetric positive-definite and therefore invertible. As a result, the RONS reduced-order equation, 
	$$ M^{(r)} \dot{\vc q} = \bxi^{(r)} + \pmb\eta^{(r)},$$
 	 is well-posed, i.e. $\dot{\vc q}$ has a unique solution. In contrast, $\Imag[M]$ is skew-symmetric and so its invertibility is not guaranteed. As a result, the reduced Lagrangian model, 
 	 $$ M^{(i)} \dot{\vc q} = \bxi^{(i)} + \pmb\eta^{(i)},$$
 	 is not necessarily well-posed.
 	 Section \ref{sec:lag_fails} contains an example where the ill-posedness of the reduced Lagrangian approach translates into its failure to produce an equation for group velocity of NLS in the stationary frame whereas RONS produces the correct group velocity. Also see Remark~\ref{rem:NLS_E} for an example where the reduced Lagrangian approach fails to produce an equation for the phase. 
\end{rem}     

Theorem~\ref{thm:nonlin_sys} shows that although RONS and the reduced Lagrangian approach are seemingly distinct, the resulting reduced-order equations
are simply the real and imaginary parts of the same complex-valued master equation. However, these real and imaginary parts do not necessarily coincide and therefore the resulting reduced-order
models can in general be different.
The following corollary gives sufficient conditions for the two reduced models, obtained by the reduced Lagrangian approach and RONS, to coincide.

\begin{cor} \label{cor:eqns_coincide}
	Let the assumptions of Theorem \ref{thm:nonlin_sys} hold. Furthermore, assume that the ansatz has the form
	\begin{equation}
	\hat{u}(x, \vc q(t)) =  \phi\bigg(\sum_{j=0}^{m} z_j(t) x^{j} \bigg),
	\label{eqn:poly_ansatz}
	\end{equation}
	where $z_j=\alpha_j + i\beta_j$, $\vc q = (\alpha_0,\cdots,\alpha_m,\beta_0,\cdots,\beta_m)$ and $\phi:\mathbb C\to \mathbb C$ is a sufficiently smooth map.
	Then the reduced-order equations produced by the reduced Lagrangian approach and RONS coincide.
\end{cor}
\begin{proof}
	See Appendix \ref{sec:thm_proof}.
\end{proof}   

Ansatz~\eqref{eqn:poly_ansatz} contains all smooth functions of polynomials in $x$ with complex coefficients.
The Gaussian ansatz~\eqref{eq:gauss_ansatz} can be put into the form~\eqref{eqn:poly_ansatz} by letting $\phi$ be the exponential function and 
defining 
$$\vc q = (\alpha_0,\alpha_1, \alpha_2, \beta_0,\beta_1,\beta_2)= (\ln(A),0,-\frac{1}{L^2}, \phi, 0, \frac{U}{L}).$$ 
Therefore, Corollary~\ref{cor:eqns_coincide} shows why the reduced Lagrangian model~\eqref{eq:param_odes}
coincides with the RONS reduced-order equation~\eqref{eq:RONS_NLS_eqn}.

\section{Stationary frame: when the reduced Lagrangian approach fails}\label{sec:lag_fails}

\begin{figure*}[t]
	\centering
	\includegraphics[width=\textwidth]{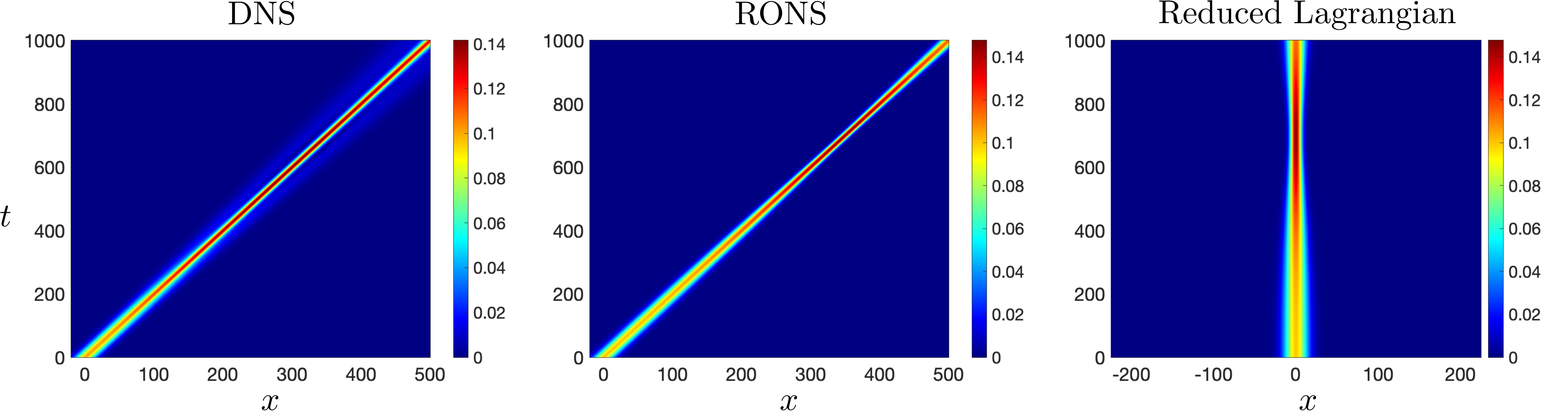}
	\caption{Comparison of the wave envelope $|u|$ obtained by DNS, RONS, and reduced Lagrangian approach. The DNS results are the full numerical simulation used as ground truth. The RONS and reduced Lagrangian results are reduced-order solutions corresponding to the translating Gaussian ansatz~\eqref{eq:gauss_ansatz_full}. Initial parameters are $A_0 = 0.1$, $L_0 = 15$, $U_0 = \phi_0 = 0$, and $x_c(0)=0$.}
	\label{fig:NLS_speeds}
\end{figure*}

We now consider the NLS equation~\eqref{eq:NLS_advection} in the stationary frame with the corresponding Lagrangian~\eqref{eq:NLS_Lag_general} with $\gamma = 1$.
We show that, although Theorem~\ref{thm:nonlin_sys} is still applicable, the appropriate ansatz does not take the form~\eqref{eqn:poly_ansatz}. As a result, Corollary~\ref{cor:eqns_coincide} is not applicable and the reduced-order equations obtained from RONS 
and the reduced Lagrangian approach do not coincide. Furthermore, the reduced Lagrangian approach
fails to predict the correct wave speed in this case, whereas RONS successfully predicts the group velocity $1/2$.

The advection term in  the NLS equation~\eqref{eq:NLS_advection} causes wave groups to translate in space. 
To account for this translation, a natural choice of ansatz is to modify the Gaussian ansatz~\eqref{eq:gauss_ansatz} so that the the center of the wave group is allowed to move. Therefore, we use the ansatz
\begin{equation}
\hat u(x, \vc q(t)) = A\exp\bigg[- \frac{(x - x_c)^2}{L^2} + i\frac{ (x - x_c)^2 U }{ L } + i \phi \bigg],
\label{eq:gauss_ansatz_advection}
\end{equation}
where the new parameter $x_c(t)$ determines the center of the wave packet. Note that, because of the translating center $x_c$, this ansatz does not admit the form of Eq.~\eqref{eqn:poly_ansatz}
required by Corollary~\ref{cor:eqns_coincide}. As a result, there is no guarantee that the reduced Lagrangian model and the RONS model would coincide.

In the case of NLS in the stationary frame, we expect the predicted dynamics of the reduced-order model to remain the same as in the co-moving frame with the exception that the wave group should translate at the constant speed $1/2$. Thus, we expect that the ODEs obtained for the parameters are the same as equation~\eqref{eq:param_odes} obtained in the stationary frame together with the additional equation  $\dot{x}_c = 1/2$.

After substituting the translating Gaussian ansatz~\eqref{eq:gauss_ansatz_advection} into the Lagrangian for NLS in the stationary frame and integrating we obtain the reduced Lagrangian
\begin{equation}
\begin{split}
\hat \Lag (\vc q, \dot{  \vc q } ) = &\frac{\sqrt{\pi } A^2 }{16L} \bigg[ \sqrt{2} \Big( -2 L^3 \dot{U}   \\ 
&+ L^2 \Big(2 \dot{L} U+U^2-8 \dot{\phi } \Big) +1\Big) -2 A^2 L^2 \bigg].
\end{split}
\label{eq:reduced_lag_advection}
\end{equation} 
which is exactly the same reduced Lagrangian~\eqref{eq:reduced_lag_comoving} obtained in the co-moving frame. 
The parameters for the wave center $x_c(t)$ and group velocity $\dot x_c(t)$ do not appear in the reduced Lagrangian~\eqref{eq:reduced_lag_advection} and so the reduced Lagrangian approach fails to produce an equation for the wave speed. We refer to Remark~\ref{rem:HPD} as to why the reduced Lagrangian approach sometimes fails to produce equations for all parameters of the ansatz.

On the other hand, RONS gives the correct results in the stationary frame.
Taking the translating Gaussian ansatz~\eqref{eq:gauss_ansatz_advection} and applying RONS to NLS in the fixed frame we obtain the parameter ODEs
\begin{align}
&\dot{A} = \frac{A U}{4 L},&  &\dot{L} = -\frac{U}{2},& \quad \dot{U} = \frac{\sqrt{2} A^2 L^2-2}{4 L^3} \nonumber \\
 &\dot{\phi} = \frac{1}{4 L^2}-\frac{5 A^2}{8 \sqrt{2}},&   &\dot{x}_c = \frac{1}{2}.&
\end{align}
These are the dynamics we expect where the wave speed is predicted to be  $1/2$ and ODEs for the other parameters remain unchanged from those in equation~\eqref{eq:param_odes} corresponding to  the co-moving frame. 

In figure \ref{fig:NLS_speeds}, we compare the true wave envelope $|u|$ obtained from direct numerical simulations (DNS) against the reduced-order solutions obtained from RONS, and reduced Lagrangian approach. 
RONS closely matches the DNS solution while the reduced Lagrangian approach fails to capture the translating wave.
We use DNS of NLS as the ground truth in order to assess accuracy of the reduced-order models.
Following Refs. \cite{cousins15,cousins16}, the DNS results are obtained using a pseudospectral methods in space and an exponential time-differencing scheme~\cite{cox02}.
We use periodic boundary conditions with a large domain size of $256 \pi$  to ensure that boundary effects do not substantially alter the results.
For all runs we use $2^{10}$ Fourier modes and a time step of $0.025$.

{\cb We also note that, as an alternative to the Gaussian ansatz, we could have chosen a hyperbolic secant ansatz proposed by Desaix et al. \cite{Desaix91},
	\begin{align}
		\hat{u} (x, \vc q (t) ) &= A\,\sech \bigg(  \frac{ x - x_c }{ L } \bigg) e^{i( x - x_c )^2 U}
	\label{eq:sech_ansatz}
	\end{align}
	with the time-dependent parameters $\vc q(t) =(A, L, U, x_c)^\top$, where $A = A^{(r)} + i A^{(i)}$ is a complex amplitude with real and imaginary parts $A^{(r)}$ and $A^{(i)}$, respectively. The remaining parameters are real-valued with $L$ being the length scale, $U$ the phase velocity, and $x_c$ the wave center.
	With the initial data $(|A(0)|^2$, $L(0)$, $U(0)$, $x_c(0)$  )$^\top$ $= (1/\sqrt{2}L_0$, $L_0$, $0$, $0  )^\top $, ansatz~\eqref{eq:sech_ansatz} is an exact soliton solution for NLS \cite{cousins15}. Other initial conditions lead to focusing or defocusing waves based on the choice of initial amplitude and length scale.
	Ideally, our reduced-order models should be able to capture the hyperbolic secant soliton as well as focusing or defocusing waves.
	
	The ansatz~\eqref{eq:sech_ansatz} does not fit the form required to apply Corollary~\ref{cor:eqns_coincide}. As a result, reduced Lagrangian approach and RONS return different reduced-order models. The resulting RONS model correctly captures the soliton solutions and in addition correctly predicts whether focusing or defocusing will occur (not shown here for brevity). 
	
	On the other hand, the reduced Lagrangian approach fails to give an equation for the wave speed $\dot{x}_c$ and, in addition, incorrectly predicts only constant amplitude waves, i.e., $|A(t)| = |A_0|$ regardless of the initial condition. More precisely, the reduced Lagrangian approach applied to the hyperbolic secant ansatz fails to predict focusing and defocusing behavior of the waves.
}

\section{Modified NLS: when the reduced Lagrangian method is inapplicable} 

\begin{figure*}[t]
	\centering
	\includegraphics[width=\textwidth]{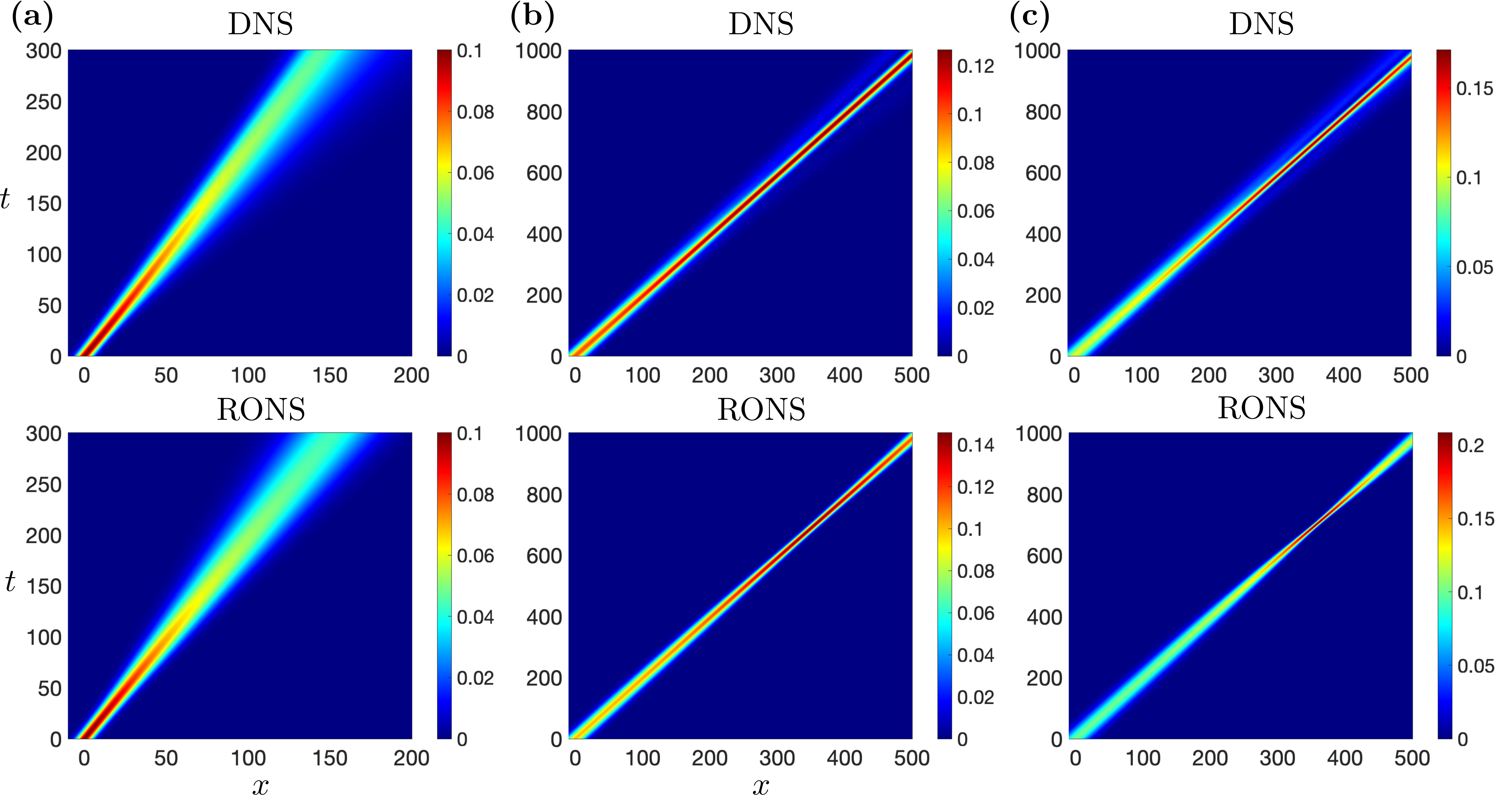}
	\caption{Comparison of evolution of $|u|$ for DNS and RONS results using the full Gaussian ansatz~\eqref{eq:gauss_ansatz_full}. Initial parameters are $U_0 = V_0 = \phi_0 = 0$ and (a) $A_0 = 0.1$, $L_0 = 5$, (b) $A_0 = 0.1$, $L_0 = 15$, (c) $A_0 = 0.1$, $L_0 = 20$. }
	\label{fig:MNLS_speeds}
\end{figure*}

\label{sec:MNLS}
The modified nonlinear Schr{\"o}dinger (MNLS) equation, first derived by Dysthe~\cite{dysthe79}, is a higher-order approximation of the full nonlinear water wave dynamics. In dimensionless units, MNLS reads
\begin{align}
\frac{ \partial  u }{ \partial  t } = &- \frac{1}{2} \frac{ \partial  u }{ \partial  x } - \frac{i}{8} \frac{ \partial^2  u }{ \partial  x^2 } + \frac{1}{16} \frac{ \partial^3  u }{ \partial  x^3 }  - \frac{i}{2} | u|^2   u  \nonumber \\
&- \frac{3}{2} | u|^2 \frac{ \partial  u }{ \partial  x } + \frac{1}{4}  u^2 \frac{ \partial  u^* }{ \partial  x } - i u \frac{ \partial \varphi }{ \partial  x } \bigg|_{z=0} ,
\label{eq:MNLS}
\end{align}
where $\varphi$ is the velocity potential which can be calculated as $\partial_x\varphi |_{z = 0} = -\mathcal{F}^{-1} \big[ |k| \mathcal{F}[ | u|^2 ] \big]/2$ where $\mathcal{F}$ denotes the Fourier transform. 
Direct numerical simulations show that the term involving the velocity potential $\varphi$ has no significant effect on the results. Therefore, for simplicity, we neglect this term when applying RONS to MNLS. 

MNLS is particularly of interest because it more accurately reproduces the wave motion observed in lab experiments~\cite{goullet2011,lo1985}. 
Cousins and Sapsis \cite{cousins15,cousins16} show that when evolving a hyperbolic secant initial condition the amplitude of focusing waves predicted by MNLS is often lower than amplitudes predicted by NLS. 
Another important distinction between the two envelope equations is that waves evolved under MNLS have a group velocity which depends on the initial wave amplitude \cite{cousins15} rather than the constant wave speed $1/2$ for waves evolved under NLS.

MNLS does not have a known Lagrangian formulation and so the reduced Lagrangian approach is not applicable to this PDE. 
However, RONS is still applicable to MNLS using the framework described in Section~\ref{sec:RONS}.
Given that the dynamics of MNLS are more complicated than NLS, we slightly modify the translating Gaussian ansatz~\eqref{eq:gauss_ansatz_advection} to include an additional imaginary term for power of $x$ so that the ansatz becomes
\begin{align}
\hat u(x, \vc q(t)) = A\exp\bigg[ & -\frac{(x - x_c)^2}{L^2} + i\frac{ (x - x_c)^2 U }{ L } \nonumber \\
&+ i (x - x_c) V + i \phi \bigg],
\label{eq:gauss_ansatz_full}
\end{align}
where $V(t)$ is a new parameter that allows the ansatz to break symmetry in space. 
Taking this modified ansatz and applying RONS to MNLS we obtain the parameter ODEs
\begin{align}
\dot{A} &= \frac{A U (2-3 V)}{8 L},  \nonumber\\
\dot{L} &= \frac{1}{4} U (3 V-2), \nonumber\\
\dot{U} &= \frac{\sqrt{2} A^2 L^2 (7 V+2)+6 V-4}{8 L^3}, \nonumber\\
\dot{V} &= -\frac{A^2 U}{2 \sqrt{2} L}, \nonumber \\
\dot{\phi} &= \frac{1}{32 L^2} \bigg( L^2 \Big(-5 \sqrt{2} A^2 (5 V+2)+6 U^2 V \nonumber \\
&\quad \quad \quad \quad \quad+4 (V-1) V^2\Big)-6 V+8 \bigg), \nonumber\\
\dot{x}_c &= \frac{1}{16} \bigg(5 \sqrt{2} A^2+\frac{3}{L^2}+3 U^2+V (3 V-4)+8\bigg).
\end{align}
\begin{figure*}[t]
	\centering
	\includegraphics[width=0.49\textwidth]{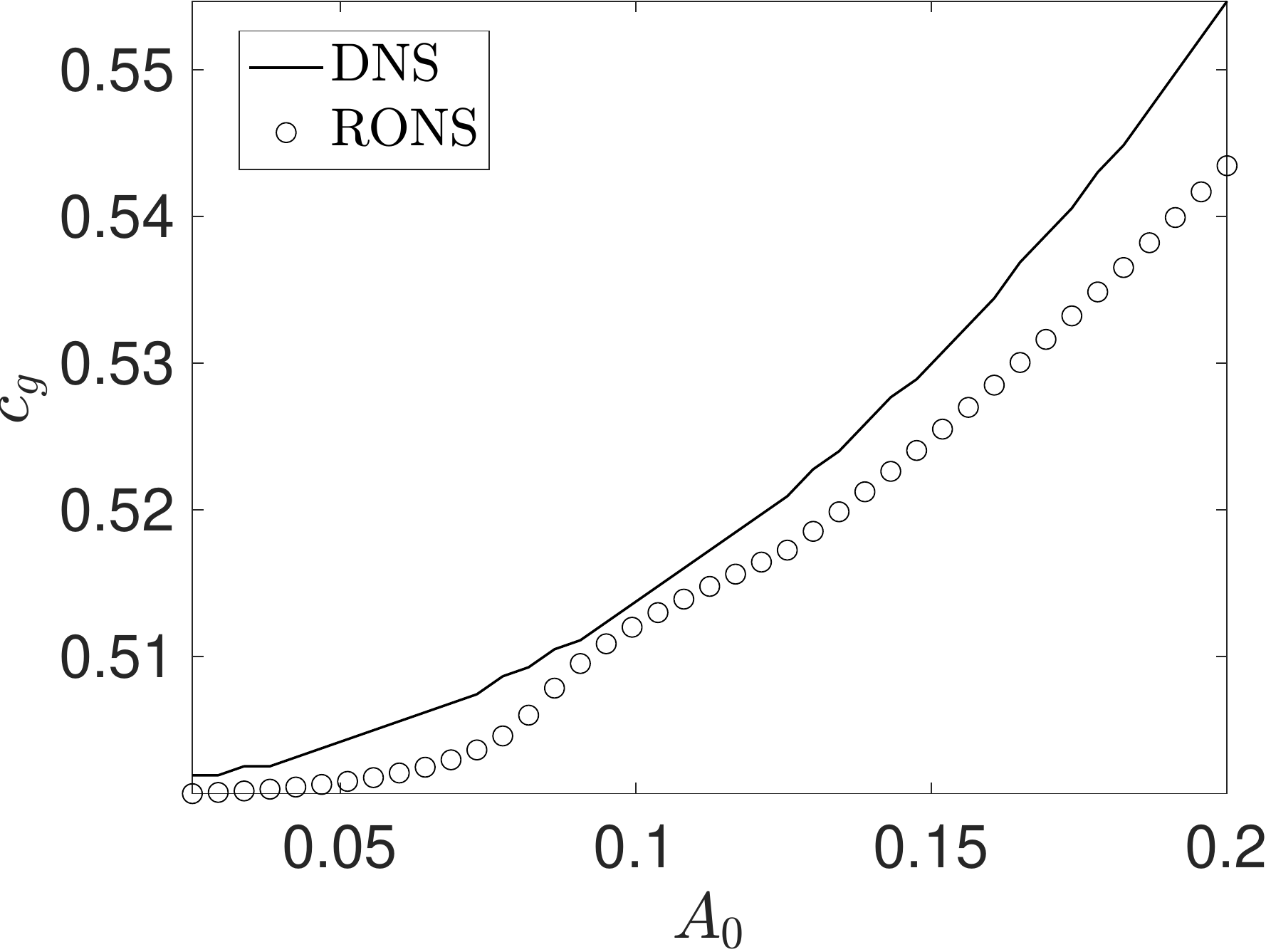}
	\includegraphics[width=0.49\textwidth]{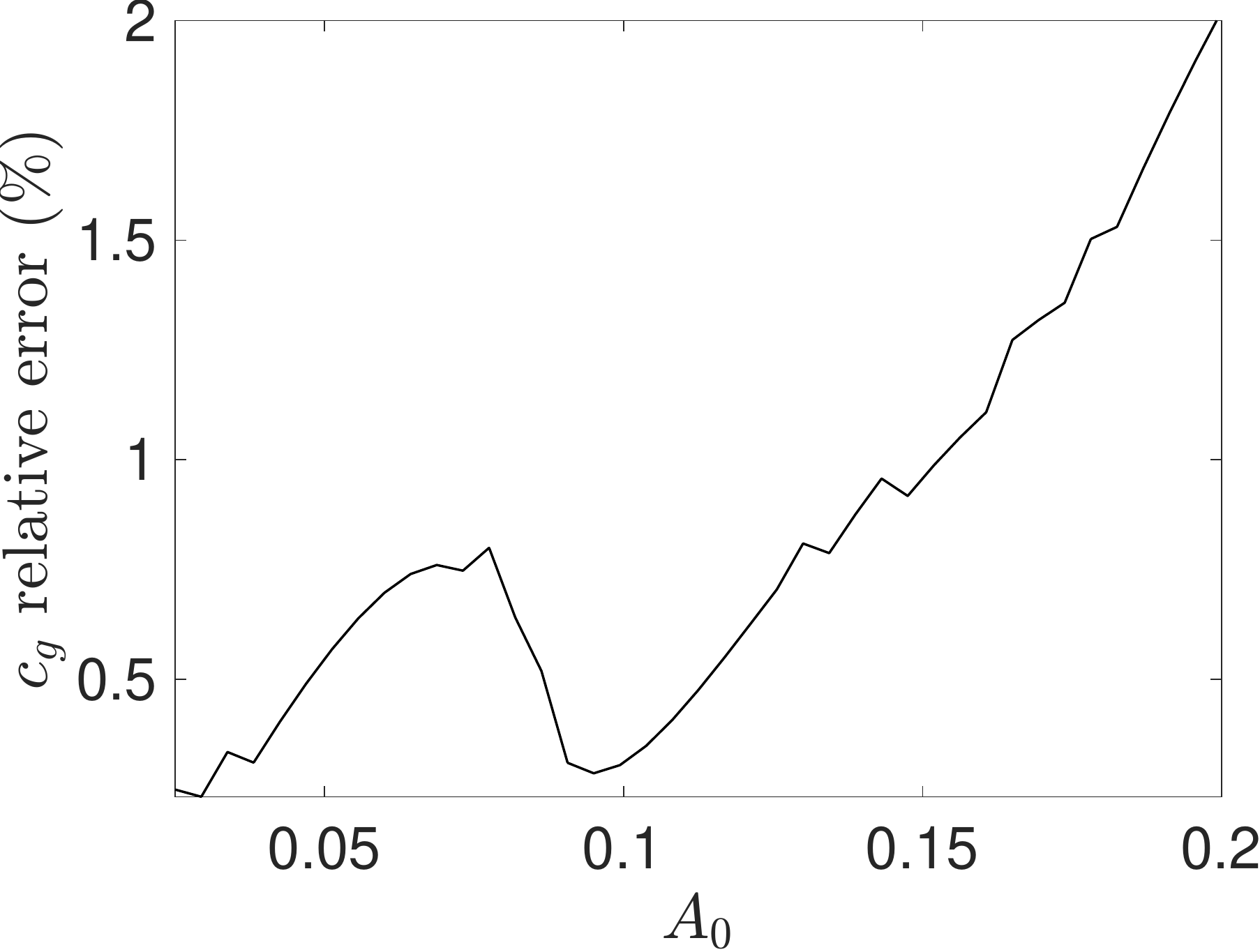}
	\caption{Comparison of group velocity $c_g$ for MNLS predicted by DNS and RONS using the full Gaussian ansatz~\eqref{eq:gauss_ansatz_full}. Initial parameters are $L_0 = 20$ and $U_0 = V_0 = \phi_0 = 0$. }
	\label{fig:MNLS_speeds_comparison}
\end{figure*}

{\cb As was the case with NLS, the reduced-order equations RONS produces for MNLS when using the modified Gaussian ansatz~\eqref{eq:gauss_ansatz_full} preserve mass without needing to explicitly enforce this conserved quantity.} 
Figure \ref{fig:MNLS_speeds} shows the evolution of the wave envelope $|u|$ for three different initial conditions. 
Comparing the RONS reduced-order solution to the DNS, we see that
in all three cases, the RONS solution accurately captures the true group velocity. Furthermore, the reduced-order solution also correctly predicts whether the wave will grow in amplitude (focusing waves) or decay (defocusing waves).
For the defocusing case the RONS reduced-order solution also does an excellent job predicting the amplitude of the wave envelope.
In the focusing case, we see that RONS provides a good approximation, although it slightly overestimates the amplitude of the wave peak and the time at which the peak takes place. 

Figure~\ref{fig:MNLS_speeds_comparison} compares the group velocity predicted by RONS against the DNS results for various initial wave amplitudes $A_0$. 
Again RONS accurately captures the group velocity with relative errors less than $2\%$. We have observed 
similar results for a wide range of waves with initial amplitudes $0<A_0\leq 0.2$ and initial length scales $5\leq L_0\leq 25$  (not shown here for brevity). In summary, RONS consistently predicts the correct group velocity for all tested initial conditions and accurately captures the behavior of defocusing waves, but slightly overestimates the maximum amplitude of focusing waves.

\section{Conclusions}
\label{sec:conclusion}
We have demonstrated the ability of RONS to develop accurate reduced-order models for NLS and MNLS equations.
For both NLS and MNLS, the RONS reduced-order models correctly predict whether initial conditions will lead to focusing or defocusing waves and approximately capture the amplitude and time of the wave peaks. Additionally, reduced solutions from RONS perfectly predict the group velocity of solutions to NLS and  accurately approximate group velocities for MNLS with at most $2\%$ relative error.

In the case of NLS, we also compared our results against reduced-order models produced from a previously developed reduced Lagrangian approach. We discovered that, although the reduced Lagrangian approach and RONS follow very different modeling philosophies, their resulting reduced-order models are intimately connected through a master equation. More precisely, the RONS model coincides with the real part of the master equation and the reduced Lagrangian model coincides with its imaginary part (see Theorem~\ref{thm:nonlin_sys}).
Furthermore, for a broad class of ansatz, the two models coincide exactly (see Corollary~\ref{cor:eqns_coincide}). 
However, RONS provides several advantages over the reduced Lagrangian approach.
We showed that the reduced Lagrangian approach was incapable of predicting wave speeds for NLS in a stationary frame whereas RONS correctly predicts a group velocity of $1/2$.
The reduced Lagrangian approach is also dependent on the PDE actually having a Lagrangian, and so it is inapplicable to many PDEs such as MNLS whereas RONS has no difficulty building a reduced-order model so long as the PDE is in the form of equation~\eqref{eq:general_pde}.
Another significant advantage of RONS is that we can easily enforce conserved quantities in solutions~\cite{anderson2021}, whereas there is no clear guarantee that obtaining reduced equations from the Lagrangian will preserve these desired quantities in solutions. These observations demonstrate the superiority of RONS for reduced-order modeling of nonlinear dispersive waves.

For future work, we would like to explore more general cases where RONS and the reduced Lagrangian approach coincide. 
We have provided sufficient conditions in Corollary~\ref{cor:eqns_coincide}, but determining necessary conditions for when the two methods coincide would likely be insightful and grant a deeper understanding of both modeling approaches.
Furthermore, the RONS framework can easily be applied to two-dimensional NLS and MNLS where we would like to examine how the reduced-order dynamics change from the one-dimensional case presented here. 
Finally, RONS can be used for reduced-order prediction of rogue waves, potentially accounting for wave-wave interactions which where neglected in similar previous studies~\cite{cousins15,cousins16,farazmand2017}.

\section*{Data availability statement}
The datasets generated during and/or analyzed during the current study are available from the
corresponding author on reasonable request.

\section*{Declarations}
Conflict of interest: The authors have no conflicts of interest to report.

\appendix

\section{Proof of Theorem \ref{thm:nonlin_sys}}
\label{sec:nonlin_sys_proof}
First we will calculate the equations for the reduced Lagrangian approach. Substituting the ansatz $\hat u(x, \vc q(t))$ into the Lagrangian~\eqref{eq:Lag_general} yields
\begin{align*}
\mathcal{L} &= \int_{\mathbb R} \bigg( \frac{1}{2} \big( i\hat u_t \hat u^\ast-i\hat u \hat u_t^\ast \big) + h(\hat u, \hat u^*, \hat u_x, \hat u_x^*) \bigg)  \id x \\
&= \int_{\R} \bigg[ \frac{1}{2} \sum_{j = 1}^{n} \bigg( i \hat{u}^* \frac{ \partial \hat{u} }{ \partial q_j }\dot{  q_j }  - i \hat{u} \frac{ \partial \hat{u}^* }{ \partial q_j }\dot{  q_j }  \bigg) +  h(\hat u, \hat u^*, \hat u_x, \hat u_x^*) \bigg] \ \id x.
\end{align*}
After a straightforward calculation and using the fact that $M$ is Hermitian, we obtain the Euler--Lagrange equations
\begin{align}
0 &= \frac{ \partial \mathcal L }{ \partial q_k }  - \frac{ \id  }{ \id  t } \frac{ \partial \mathcal L }{ \partial \dot{q_k} }  \nonumber \\
&= \int_{\R} \bigg[ \sum_{j = 1}^{n} i\bigg(  \frac{ \partial \hat{u}^* }{ \partial q_k } \frac{ \partial \hat{u} }{ \partial q_j }\dot{  q_j }  -  \frac{ \partial \hat{u} }{ \partial q_k } \frac{ \partial \hat{u}^* }{ \partial q_j }\dot{  q_j }  \bigg)  \nonumber \\
& \quad \quad \quad + \bigg( \frac{ \partial h }{ \partial \hat u } \frac{ \partial \hat u }{ \partial q_k }  + \frac{ \partial h }{ \partial \hat u^* } \frac{ \partial \hat u^* }{ \partial q_k } \bigg) \nonumber \\
&\quad \quad \quad  + \bigg( \frac{ \partial h }{ \partial \hat u_x } \frac{ \partial \hat u_x }{ \partial q_k } + \frac{ \partial h }{ \partial \hat u_x^* } \frac{ \partial \hat u_x^* }{ \partial q_k } \bigg) \bigg] \ \id x \nonumber \\
&= -2 \sum_{j = 1}^{n}\Imag[M_{jk}] \dot q_j - 2\Imag[ \xi_k] - 2\Imag[  \eta_k] \nonumber \\
&= 2\sum_{j = 1}^{n} \Imag[M_{kj}] \dot q_j - 2\Imag[ \xi_k] - 2\Imag[  \eta_k],
\label{eq:Lag_eqns}
\end{align}
for $k=1,2,\cdots,n$.
In order to calculate the equations for RONS, we recall that the PDE corresponding to Lagrangian~\eqref{eq:Lag_general} is
\begin{equation}
u_t = -i \frac{\partial }{ \partial x } \frac{\partial  h }{ \partial u_x^* } + i \frac{\partial h }{ \partial u^* }.
\label{eq:Lag_PDE}
\end{equation}
Applying RONS to~\eqref{eq:Lag_PDE} and integrating by parts yields the equations
\begin{align}
0 &=  \frac{ \partial  }{ \partial \dot{q_k} }  \| \hat u_t - F(\hat u) \|^2_{L^2} \nonumber \\
&= \int_{\R} \bigg[ \sum_{j=1}^{n} \bigg(  \frac{ \partial \hat{u}^* }{ \partial q_k } \frac{ \partial \hat{u} }{ \partial q_j }\dot{  q_j }  +   \frac{ \partial \hat{u} }{ \partial q_k } \frac{ \partial \hat{u}^* }{ \partial q_j }\dot{  q_j }  \bigg)  \nonumber \\
&\quad \quad \quad + i\bigg(   \frac{ \partial \hat u }{ \partial q_k }  \frac{\partial h }{ \partial \hat u } -  \frac{ \partial \hat u^* }{ \partial q_k } \frac{\partial h }{ \partial \hat u^* } \bigg) \nonumber \\
&\quad \quad \quad + i \bigg(  \frac{ \partial \hat u_x }{ \partial q_k } \frac{\partial  h }{ \partial \hat u_x } -  \frac{ \partial \hat u_x^* }{ \partial q_k } \frac{\partial  h }{ \partial \hat u_x^* }  \bigg)  \bigg] \id x \nonumber \\
&= 2\sum_{j = 1}^{n} \Real[M_{jk}] \dot q_j  - 2\Real[ \xi_k ] - 2\Real[ \eta_k ] \nonumber\\
&= 2\sum_{j = 1}^{n} \Real[M_{kj}] \dot q_j  - 2\Real[ \xi_k ] - 2\Real[  \eta_k ].
\label{eq:RONS_eqns}
\end{align}
Therefore, if we consider the system
\begin{equation}
M \dot{ \vc q } = \pmb \xi + \pmb \eta ,
\end{equation}
taking the imaginary part gives the reduced Lagrangian approach~\eqref{eq:Lag_eqns} and taking the real part gives RONS~\eqref{eq:RONS_eqns}.

\section{Proof of Corollary \ref{cor:eqns_coincide}} \label{sec:thm_proof}
We first note that the ansatz~\eqref{eqn:poly_ansatz} satisfies the condition
\begin{align}
i\dfrac{\partial \hat u }{ \partial \alpha_k} = \dfrac{\partial \hat u }{ \partial \beta_k},
\label{eq:derivative_cond}
\end{align}
which is a key property that will be used throughout the proof. Defining $\vc q = (\alpha_0,...,\alpha_m,\beta_0,...,\beta_m)$ we have
\begin{align*}
\frac{ \partial \hat u }{ \partial q_k } &= \begin{cases}
x^k \phi^\prime, &1 \leq k \leq m+1 \medskip \\
ix^{k-m-1} \phi^\prime, &m+1 < k \leq 2(m+1).
\end{cases}
\end{align*}
By Theorem \ref{thm:nonlin_sys} the reduced Lagrangian and RONS equations can be derived from the system of equations~\eqref{eq:nonlin_sys}.
Examining the left-hand side of the system, we write $M$ as a block matrix 
\begin{equation}
M = \begin{pmatrix}
B &\quad C^\ast\\
C &\quad  D
\end{pmatrix},
\end{equation}
where the matrices $B,C,D\in\C^{ (m+1) \times (m+1)}$ have entries given by 
\begin{align*}
	B_{lj} &= \ds \int_{\R } \dfrac{ \partial \phi }{ \partial \alpha_l }  \dfrac{ \partial \phi^* }{ \partial \alpha_j }  \ \id x, \quad C_{lj} = \ds \int_{\R } \dfrac{ \partial \phi }{ \partial \beta_l }  \dfrac{ \partial \phi^* }{ \partial \alpha_j }  \ \id x,\\
	 D_{lj} &= \ds \int_{\R } \dfrac{ \partial \phi }{ \partial \beta_l }  \dfrac{ \partial \phi^* }{ \partial \beta_j }  \ \id x.
\end{align*}
 
Matrix $B$ is clearly Hermitian, and if condition~\eqref{eq:derivative_cond} is satisfied then $B= -iC =D$. Thus, we can write $M$ as 
\begin{equation}
M = \begin{pmatrix}
B & C^*\\
C & D
\end{pmatrix} = 
\begin{pmatrix}
B & -iB\\
iB & B
\end{pmatrix} .
\end{equation}
Now we consider the vectors on the right-hand side of system~\eqref{eq:nonlin_sys}.
If condition~\eqref{eq:derivative_cond} is satisfied then $\pmb \xi$ can be written as 
\begin{align*}
\xi_k  &=  \begin{cases}
\ds -i \int_{\R } \dfrac{\partial \phi }{ \partial \alpha_k} \dfrac{\partial h}{ \partial \hat u }  \ \id x,  & 1 \leq k \leq m+1 \medskip  \medskip \\
\ds \int_{\R }  \dfrac{\partial \phi }{ \partial \alpha_{k-m-1}} \dfrac{\partial h}{ \partial \hat u }   \ \id x,  &  m+1 < k \leq 2(m+1) .
\end{cases} 
\end{align*} 
So, we may write $\pmb \xi$ as two vectors of length $m+1$ stacked on top of each other, i.e. $\pmb \xi = (\pmb \xi_1, i \pmb \xi_1)^\top$. The case for $\pmb \eta$ is similar so that we may write $\pmb \eta = (\pmb \eta_1, i \pmb \eta_1)^\top$.
Defining the vectors $\pmb \alpha := (\alpha_0,...,\alpha_m)^\top$ and  $\pmb \beta := (\beta_0,...,\beta_m)^\top$ we write the nonlinear system for the parameter evolution~\eqref{eq:nonlin_sys} as 
\begin{equation}
\begin{pmatrix}
B & -iB\\
iB & B
\end{pmatrix}\begin{pmatrix}
\dot{ \pmb \alpha } \\
\dot{ \pmb \beta }
\end{pmatrix}  =  \begin{pmatrix}
\pmb \xi_1 \\
i\pmb \xi_1
\end{pmatrix} + \begin{pmatrix}
\pmb \eta_1 \\
i\pmb \eta_1
\end{pmatrix} .
\end{equation}
Taking the real part gives RONS which yields
\begin{align}\label{eq:Re}
\underbrace{  \begin{pmatrix}
	\Real[B] & \Imag[ B ]\\
	-\Imag[ B ] & \Real[ B ]
	\end{pmatrix} }_{ M^{ (r) } } \underbrace{ \begin{pmatrix}
	\dot{ \pmb \alpha } \\
	\dot{ \pmb \beta }
	\end{pmatrix}}_{ \dot{  \vc q } }   &=   \underbrace{  \begin{pmatrix}
	\Real[ \pmb \xi_1 ] \\
	-\Imag [ \pmb \xi_1 ] 
	\end{pmatrix}}_{ \pmb \xi^{ (r) } }  +  \underbrace{  \begin{pmatrix}
	\Real[ \pmb \eta_1 ] \\
	-\Imag[ \pmb \eta_1 ]
	\end{pmatrix}}_{ \pmb \eta^{ (r) } }  
\end{align}
Taking the imaginary part gives the Lagrangian approach which yields
\begin{align} \label{eq:Im}
\underbrace{ \begin{pmatrix}
	\Imag[ B ] & -\Real[ B ]\\
	\Real[ B ]  &  \Imag[ B ]
	\end{pmatrix} }_{ M^{ (i) } } \underbrace{ \begin{pmatrix}
	\dot{ \pmb \alpha } \\
	\dot{ \pmb \beta }
	\end{pmatrix}}_{ \dot{  \vc q } }    &=  \underbrace{ \begin{pmatrix}
	\Imag[ \pmb \xi_1 ]  \\
	\Real[ \pmb \xi_1 ]
	\end{pmatrix}}_{ \pmb \xi^{ (i) } }   + \underbrace{ \begin{pmatrix}
	\Imag[ \pmb \eta_1] \\
	\Real[ \pmb \eta_1]
	\end{pmatrix}}_{ \pmb \eta^{ (i) } }.
\end{align}
Note that equations~\eqref{eq:Re} and~\eqref{eq:Im} are identical since both lead to the set of equations
\begin{equation}
\begin{cases}
\Real[B] \dot{ \pmb \alpha } + \Imag[B] \dot{ \pmb \beta } = \Real[\pmb\xi_1] + \Real[\pmb \eta_1]\\
\Imag[B] \dot{ \pmb \alpha } - \Real[B] \dot{ \pmb \beta } = \Imag[\pmb\xi_1] + \Imag[\pmb \eta_1]
\end{cases}
\end{equation}
Therefore the reduced-order equations obtains through RONS coincide with those obtained by the reduced Lagrangian method.


%
%


\end{document}